\documentclass[12pt]{amsart}
\usepackage{graphicx, amsmath, hyperref, enumitem} 
\usepackage[svgnames]{xcolor}
\usepackage{csquotes}
\usepackage{appendix}

\usepackage[
backend=biber,
style=numeric,
sorting=nyt,
]{biblatex}
\newlength{\minipagewidth}
\makeatletter
\@addtoreset{equation}{section}
\renewcommand{\theequation}{\thesection.\@arabic\c@equation}
\newcommand\Nopagebreak{\@nobreaktrue\nopagebreak}
\makeatother

\newtheorem*{thm*}{Theorem}
\newtheorem*{conj*}{Conjecture}

\theoremstyle{definition}

\newtheorem*{axiom}{Axiom}

\newtheorem*{example}{Example}

\newtheorem*{exs}{Examples}
\newtheorem*{defn*}{Definition}

\usepackage[framemethod=TikZ]{mdframed}
\mdfdefinestyle{theorem}{%
    linecolor=black,
    outerlinewidth=2pt,
    roundcorner=8pt,
    innertopmargin=-2pt,
    innerbottommargin=6pt,
    innerrightmargin=5pt,
    innerleftmargin=5pt,
    leftmargin=-7pt,
    rightmargin=-7pt,
    backgroundcolor=yellow!15!white}
\mdfdefinestyle{theoremm}{%
    linecolor=black,
    outerlinewidth=2pt,
    roundcorner=8pt,
    innertopmargin=6pt,
    innerbottommargin=6pt,
    innerrightmargin=5pt,
    innerleftmargin=5pt,
    leftmargin=-7pt,
    rightmargin=-7pt,
    backgroundcolor=yellow!15!white}
\mdfdefinestyle{conjecture}{%
    linecolor=black,
    outerlinewidth=2pt,
    roundcorner=8pt,
    innertopmargin=-2pt,
    innerbottommargin=6pt,
    innerrightmargin=5pt,
    innerleftmargin=5pt,
    leftmargin=-7pt,
    rightmargin=-7pt,
    backgroundcolor=red!15!white}
\mdfdefinestyle{definition}{%
    linecolor=black,
    outerlinewidth=2pt,
    roundcorner=8pt,
    innertopmargin=-2pt,
    innerbottommargin=6pt,
    innerrightmargin=5pt,
    innerleftmargin=5pt,
    leftmargin=-7pt,
    rightmargin=-7pt,
    backgroundcolor=blue!10!white}
\mdfdefinestyle{proof}{%
    linecolor=black,
    outerlinewidth=2pt,
    roundcorner=8pt,
    innertopmargin=6pt,
    innerbottommargin=6pt,
    innerrightmargin=5pt,
    innerleftmargin=5pt,
    leftmargin=-7pt,
    rightmargin=-7pt,
    backgroundcolor=green!15!white}
 
\mdfdefinestyle{example}{%
    linecolor=black,
    outerlinewidth=2pt,
    roundcorner=8pt,
    innertopmargin=-2pt,
    innerbottommargin=6pt,
    innerrightmargin=5pt,
    innerleftmargin=5pt,
    leftmargin=-7pt,
    rightmargin=-7pt,
    backgroundcolor=orange!15!white}

    \mdfdefinestyle{nonexample}{%
    linecolor=black,
    outerlinewidth=2pt,
    roundcorner=8pt,
    innertopmargin=-2pt,
    innerbottommargin=6pt,
    innerrightmargin=5pt,
    innerleftmargin=5pt,
    leftmargin=-7pt,
    rightmargin=-7pt,
    backgroundcolor=yellow!15!white}

\mdfdefinestyle{aside}{%
    linecolor=black,
    outerlinewidth=2pt,
    roundcorner=8pt,
    innertopmargin=-2pt,
    innerbottommargin=6pt,
    innerrightmargin=5pt,
    innerleftmargin=5pt,
    leftmargin=-7pt,
    rightmargin=-7pt,
    backgroundcolor=white}
    
\surroundwithmdframed[style=theorem]{thm}
\surroundwithmdframed[style=theorem]{thm*}
\surroundwithmdframed[style=theoremm]{thmm}
\surroundwithmdframed[style=theorem]{cor}
\surroundwithmdframed[style=theorem]{lemm}
\surroundwithmdframed[style=theorem]{prop}
\surroundwithmdframed[style=conjecture]{conj}
\surroundwithmdframed[style=proof]{proof}
\surroundwithmdframed[style=definition]{defn}
\surroundwithmdframed[style=definition]{defn*}
\surroundwithmdframed[style=definition]{axiom}
\surroundwithmdframed[style=example]{example}
\surroundwithmdframed[style=nonexample]{nonexample}
\surroundwithmdframed[style=example]{examples}
\surroundwithmdframed[style=example]{exs}
\surroundwithmdframed[style=aside]{aside}

\def\Z{\mathbb{Z}}
\def\N{\mathbb{N}}

\title{Reflective Groupwork For Introductory Proof-Writing Courses}
\author{Jennifer Pi}
\author{Christopher Davis}
\author{Yasmeen Baki}
\author{Alessandra Pantano}
\address{University of California, Irvine}
\email{daviscj@uci.edu, ybaki@uci.edu, apantano@uci.edu, jspi@uci.edu}

\begin{document}

\begin{abstract}
    We discuss two proof evaluation activities meant to promote the acquisition of learning behaviors of professional mathematics within an introductory undergraduate proof-writing course. These learning behaviors include the ability to read and discuss mathematics critically, reach a consensus on correctness and clarity as a group, and verbalize what qualities ``good" proofs possess. 
    The first of these two activities involves peer review and the second focuses on evaluating the quality of internet resources.
    All of the activities involve groupwork and reflective discussion questions to develop students' experience with these practices of professional mathematics.
\end{abstract}

\maketitle

\begin{displayquote}
    \textit{The understanding of learning and teaching mathematics...support[s] a model of participating in a culture rather than a model of transmitting knowledge...the core effects as emerging from the participation in the culture of a mathematics classroom will appear on the metalevel mainly and are ``learned" indirectly.} 
    \begin{flushright}
    - Heinrich Bauersfeld \cite{bauersfeld1993teachers}
    \end{flushright}
\end{displayquote}

\section{Introduction}
As preparation for rigorous proof-based mathematics courses such as abstract algebra and real analysis, many universities now offer an ``Introduction to Mathematical Reasoning", or transition-to-proofs, course. 
The purpose of such a course is to provide students with the tools to read and understand abstract mathematics, learn formal mathematical language and definitions, and of course,  gain the ability to write proofs of mathematical statements.
While the focus of the course is on proof-writing, there are ancillary goals of introducing students to the  practices of professional mathematics. Sociomathematical norms, defined by Yackel and Cobb \cite{sociomathematicalnorms} as ``normative aspects of mathematical discussions that are specific to students' mathematical activity'', are implicitly established throughout students' mathematical experiences. However, professional mathematics practices are not typically explicitly emphasized in courses, and the introduction to proofs course is usually the first exposure most students have to these professional practices.

In this paper, we present two proof evaluation activities intended to promote student reflection on some of the norms of professional mathematics.
The first is a peer review activity emphasizing the traits that a ``good" proof possesses (this is an intentionally vague prompt), and fits into a long history of peer review in teaching mathematics; see for example \cite{peerfeedback, learningfrommistakes2011}.
The second activity focuses on evaluating online resources and reflecting on what qualities make a proof easier to understand.
Both activities emphasize the collaborative nature of mathematical work via required group submission.
Having students practice collaboration with their classmates is not only an added instrument in their mathematical toolbox, but also helps to build relationships with their peers and thus retain motivation for the course. The hope is that such groupwork builds a stronger sense of community amongst the students, and promotes the emergence of sociomathematical norms which align with professional mathematical practices via the active practice of \emph{doing} (as opposed to learning) mathematics collaboratively.

\section{Proof Evaluation 1: What Makes a ``Good" Proof?}


\subsection{Description of Activity}
The first activity focuses on peer review of proofs. 
The main learning goals of this activity are for students to
\begin{enumerate}
    \item \label{SLO correctness and clarity} Read mathematics critically, for both correctness and clarity.
    \item Discuss and critique peer work.
    \item Reach a consensus of what makes a proof ``good".
\end{enumerate}

For the activity, each student is first asked to prove or disprove the following:
\begin{displayquote} 
    \emph{(Prompt)}     If $x$ is a prime number, then $\sqrt{x}$ is irrational.
\end{displayquote}

 Then students are instructed to: 
\begin{enumerate}[label=(\alph*), ref=(\alph*)]
\item Form groups of three (or four) and exchange proofs with their group members;
\item Perform a peer review of their group members' proofs, checking both for correctness and some kind of clarity or `aha' feeling;

\item Meet with their group members to answer some discussion questions, and rewrite one of the individual proofs that the group agreed had some confusing or incorrect aspects, then reflect upon which qualities make a proof ``good";

\item Complete a group submission including the individual proofs originally written by the group members, along with the answers to the group discussion questions.
\end{enumerate}

This activity was given as a take-home assignment; students had a week to complete the activity. Instructions were given through the course management platform (Canvas), though the assignment was also announced during several in-person lectures. Submissions were also online via Canvas. This activity was assigned during the third week of a 10-week quarter, when students had just finished covering basic set theory, propositional functions and quantifiers, and the three main methods of proof (direct, contradiction, contrapositive). 
To view a copy of the activity as given to the students, please see Appendix ~\ref{appendix: Proof Eval 1}.


\subsection{A Success Story: Rewriting a Proof as a Group} \label{sec: PE1 success story}
We discuss here a submission completed by a group of four students, who herein will be called Nora, David, Yonathan and Oliver (all names are pseudonyms). 
The first part of the submission contains every student's initial proof, as per instructions. 
We note that all students attempted to do a proof by contradiction: they assumed that $x$ is a prime number and that $\sqrt{x}$ is rational (say $\sqrt{x} = \frac{a}{b}$, with $a$ and $b$ co-prime), and showed a contradiction by proving that the fraction $\frac{a}{b}$ cannot be in reduced terms. While all students chose the same approach, the quality of their proofs varied, both in regards to logic and style. 
Bridging these differences in logic and style as a group and reaching a consensus on the clarity of each proof was the objective of the activity.
We include below images of the four original proofs. 
\\

\begin{figure}[htp]
    \centering
\includegraphics[width=11cm]{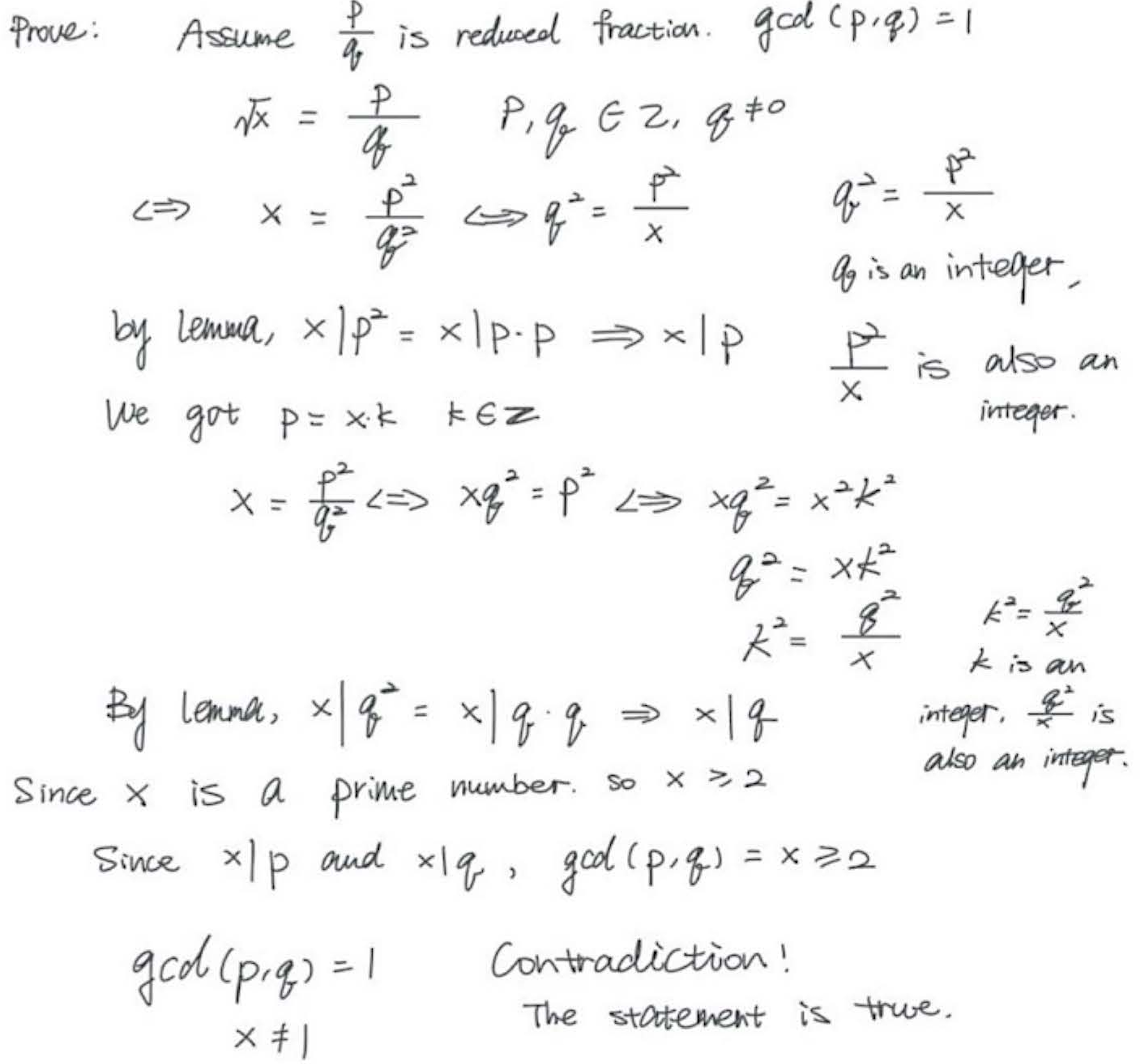}
    \caption{David's proof}
    \label{fig:david}
     \end{figure}  
    
    \begin{figure}[htp]
    \centering
     \end{figure}
     \begin{figure}[htp]
    \centering
\includegraphics[width=11cm]{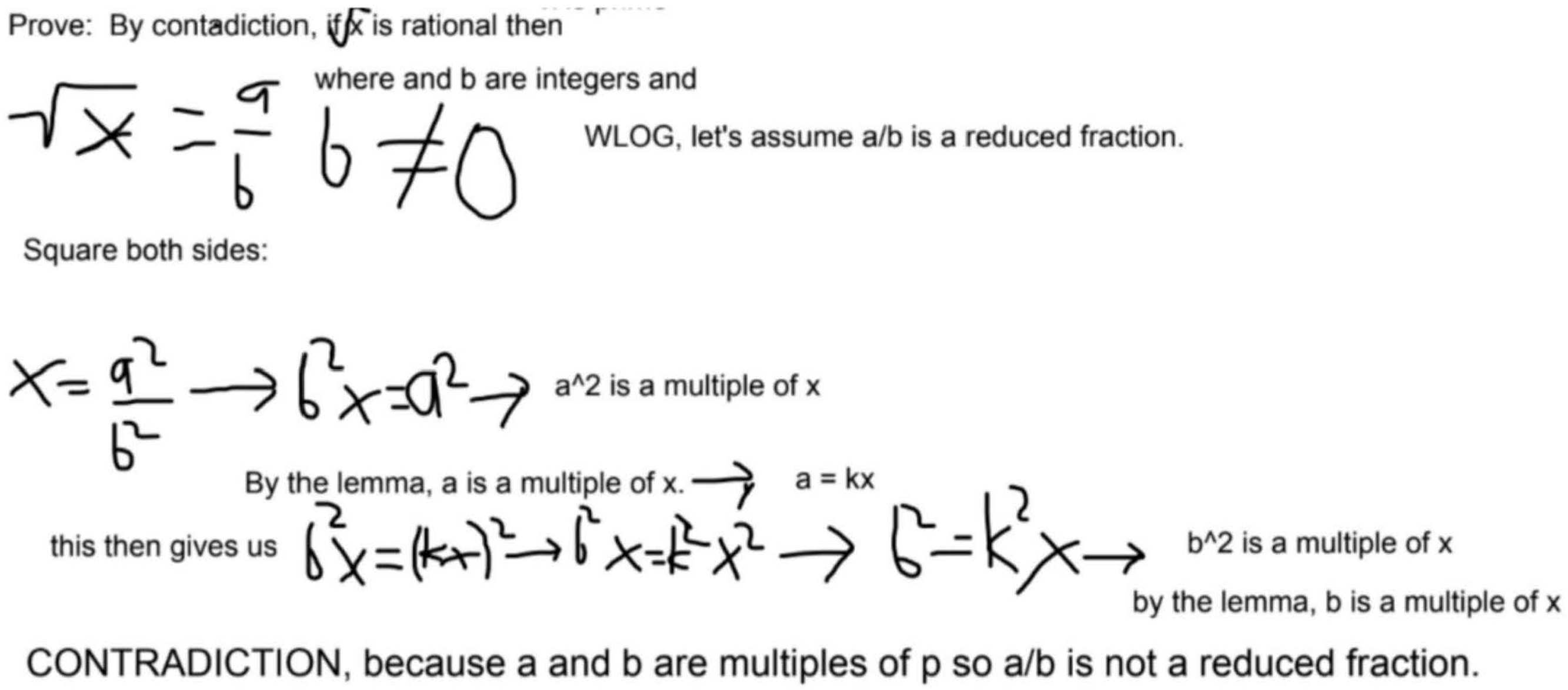}
    \caption{Oliver's proof}
    \label{fig:oliver}
     \end{figure}  
     
\begin{figure}[htp]
    \centering
\includegraphics[width=11cm]{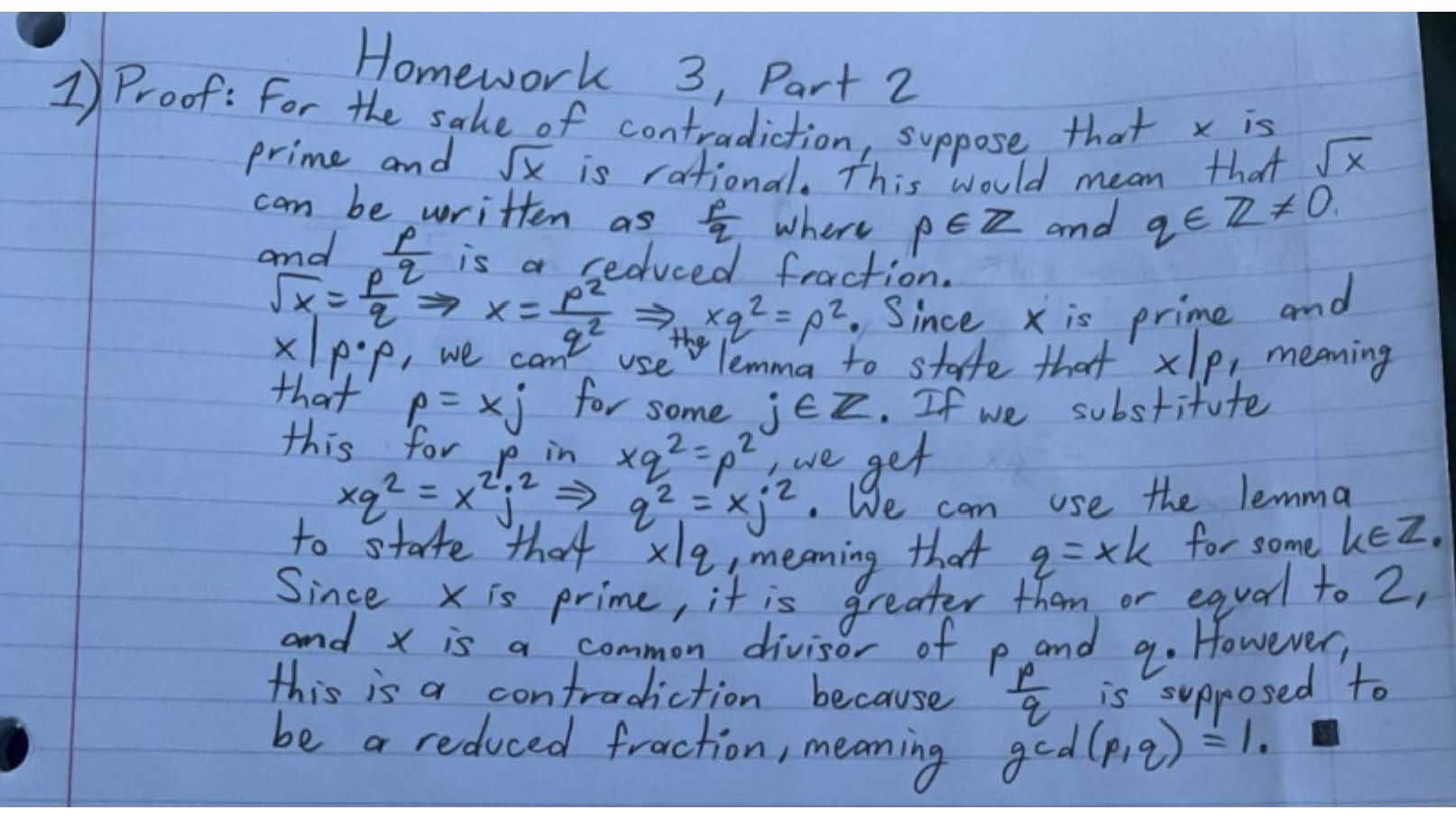}
    \caption{Yonathan's proof}
    \label{fig:yonathan}
\end{figure}

     \begin{figure}[htp]
\centering\includegraphics[width=11cm]{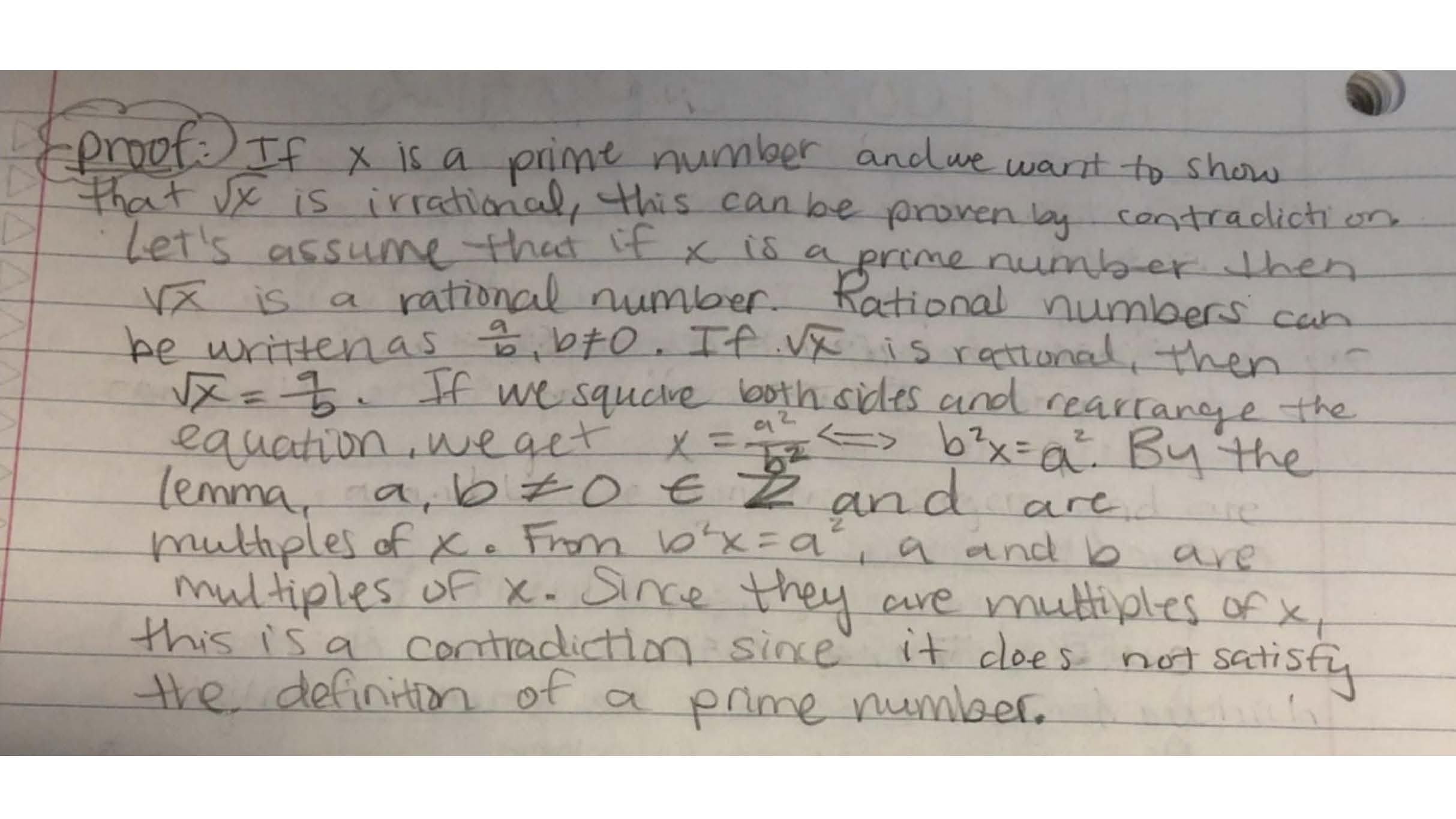}
    \caption{Nora's proof}
    \label{fig:nora}
    \end{figure}

\begin{figure}[htp]
    \centering
\includegraphics[width=10cm]{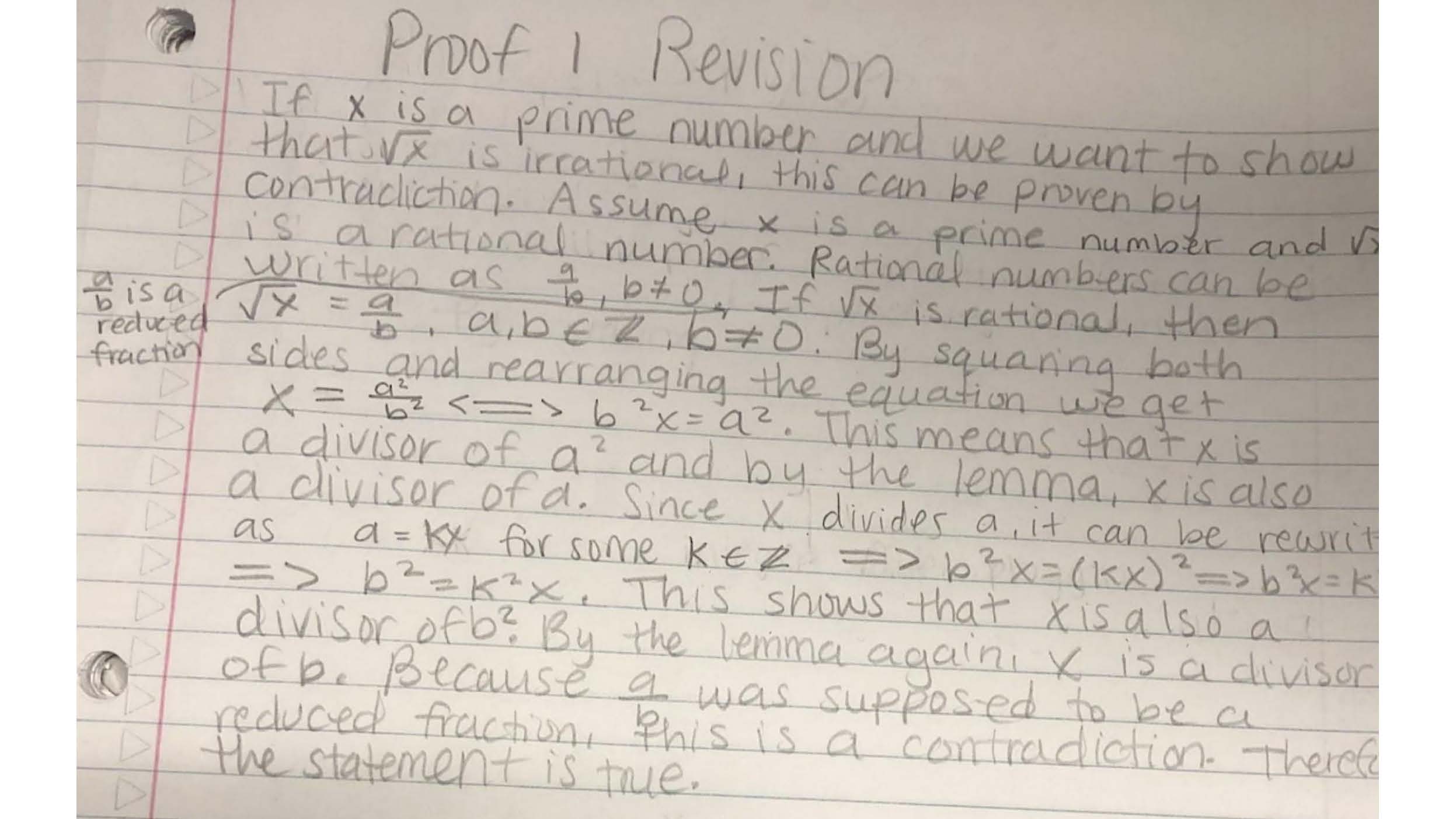}
    \caption{Group revision}
    \label{fig:revision}
\end{figure}

After presenting the four original proofs, the group was asked to
\begin{enumerate}[label=(\alph*), ref=(\alph*)]
    \item identify a part of one of the  proofs that was confusing; 
    \item explain (as a group) why some members were confused by this portion of the selected proof; 
    \item provide a better explanation of that portion of the proof as a group.
\end{enumerate}
Nora, David, Yonathan and Oliver selected Nora's proof, stating that it may benefit from additional clarity and justification. The portion of the proof they found most confusing was the explanation of why  the rational number $\sqrt{x}$ is not a reduced fraction. Quoting the students: 
\begin{displayquote}
    \emph{(Groupwork Quote)}  Our group found that proof 1 was the most confusing. It seems as if there needed to be more justifications and explanations, like including that $p/q$ is a reduced fraction. Some of my sentences were not clear, such as my opening.
\end{displayquote}
Then, they included the following write-up as the group's attempt to clarify the portion of Nora's proof which they found confusing (see Figure \ref{fig:revision}).

We note with satisfaction that the group implemented three major improvements in the revision: 
\begin{enumerate}
\item Although both the original proof by Nora and the group revision argued by contradiction, the formulation of the hypothesis in Nora's original proof was incorrect: 
\vskip0.2cm
\begin{displayquote}
    \emph{(Quote from Nora)} Let's assume that if $x$ is a prime number then $\sqrt{x}$ is a rational number.
\end{displayquote}
\vskip.2cm 

    \noindent
In the group revision, the error was fixed. The new proof correctly states:
\vskip.2cm 
\begin{displayquote}
    \emph{(Group Revision Quote)} Let's assume that $x$ is a prime number and $\sqrt{x}$ is a rational number. 
\end{displayquote}

    \item The original proof by Nora lacked argumentation as to why $x$ had to divide both $a$ and $b$. Quoting Nora's work:  
    \vskip.2cm 
\begin{displayquote}
    \emph{(Quote from Nora)} By the lemma, $a,\,b\neq 0 \in \mathbb{Z}$ and $a$, $b$ are multiples of $x$. From $b^2x=a^2$, $a$ and $b$ are multiples of $x$.
\end{displayquote}
\vskip.2cm 

    \noindent
The group revision includes a much more thorough explanation of why $x$ divides both $a$ and $b$:
\vskip.2cm 
\begin{displayquote}
    \emph{(Group Revision Quote)}
    This [$b^2x=a^2$] means that $x$ is a divisor of $a^2$ and by the lemma, $x$ is also a divisor of $a$. Since $x$ divides $a$, it can be rewritten as $a=kx$ for some $k\in \Z$.  $\Rightarrow \, b^2x= (kx)^2 \,
    \Rightarrow \, b^2x= k^2x^2 \, \Rightarrow \,  b^2= k^2x$.  This shows that $x$ is also a divisor of $b^2$. By the lemma again, $x$ is also a divisor of $b$. \end{displayquote}
\vskip.2cm 

    \item In her original proof, Nora claimed that the fact that $x$ divides both $a$ and $b$ contradicts the definition of prime: 
\vskip.2cm 
\begin{displayquote}
    \emph{(Quote from Nora)} Since they are multiple of $x$, this is a contradiction since it does not satisfy the definition of a prime number.
\end{displayquote}
\vskip.2cm 

\noindent 
The group revision, on the other hand, correctly attributes the contradiction to the fact that the (reduced) fraction $\frac{a}{b}$ is not in lowest terms: 
\vskip.2cm 
\begin{displayquote}
    \emph{(Group Revision Quote)}  Because $\frac{a}{b}$ was supposed to be a reduced fraction, this is a contradiction.
\end{displayquote}
    \vskip.2cm 
    
\end{enumerate}

Other facets of the group revision which improve upon Nora's original work are the observation that $b$ is nonzero and the clear conclusion (``Therefore the statement is true.''). \\

The second part of the group submission explicitly asked students to consider what qualities a ``good" proof inherently has:
\begin{displayquote}
    \emph{(Prompt)}  Write a few sentences about what your group agrees upon as qualities of a ``good'' proof, and give examples of these qualities from the proofs on the previous page.
\end{displayquote}


Stepping away from Nora, David, Yonathan and Oliver's group submission, we report the  general attributes of a ``good'' proof mentioned collectively by the students in our two classes, organized by category (and sub-category). The number next to a sub-category indicates its frequency.  
\begin{itemize}
\item \emph{Logical contents of the proof}
\begin{itemize}
\item Clear logical order (9)
\item No gaps in the logic or cases omitted (level of rigor) (3)
\end{itemize}

\item \emph{How the proof is written: narrative structure}
\begin{itemize}
\item Restating goal at end of proof/conclusion (5)
\item Recalling definitions (3) - Interestingly, one submission specifically mentions that it is better to \emph{not} recall definitions like ``prime number''
\item Preview of how it will be proved (3)
\item Minimize calculations (1)
\end{itemize}
\item \emph{How the proof is written: terminology}
\begin{itemize}
\item Defining variables (4)
\item Use of complete sentences (2)
\item Useful symbols (1)
\item Correct notation (1) \end{itemize}
\item \emph{How the proof is written: amount of explanation}
\begin{itemize}
\item Clear definitions/axioms/assumptions (8)
\item Tidy and concise (6)
\item Not overly complicated (2)
\item Use neither too much nor too little explanation (1)
\item Minimum amount of words necessary (1)
\end{itemize}
\end{itemize} 

In summary, our students agree that a ``good'' proof  is tidy and concise. It starts by stating definitions/axioms/assumptions, follows a clear logical order, and ends by restating the goal at end of the proof. Defining variables and recalling definitions are generally regarded as positive practices in writing a proof. 
It is interesting to note that there is some tension between student responses: some groups noted that recalling too many definitions that should be ``well-known" can conflict with the quality of being \emph{tidy and concise}, while others wished for additional clarity via clear definitions and a preview of the proof before the actual logical content. 
Any professional mathematician recognizes that this type of benchmark for how much to write for a ``good" proof shifts with both the author and the audience.

\subsection{Implementation and Adaptation} 
In this section, we share some shortcomings of the current version of this activity, as presented in Appendix \ref{appendix: Proof Eval 1}, and we propose some suggestions for future implementations.

\begin{enumerate} \item \emph{One of the questions in Proof Evaluation 1 is grounded in a deficit perspective.} 

We feel that the biggest shortcoming of this activity is its potential to make 
students feel discouraged due to the critiques of their proofs. Upon presenting the group members' original proofs, the students are asked to identify a part of one of the proofs that was confusing.
Having a portion of your original proof be selected as the most confusing one from the group may be embarrassing. Similarly, seeing your peers come together to ``fix'' your proof may be discomfiting. We already saw an example of these feelings emerge in  Nora, David, Yonathan and Oliver's submission, when Nora wrote as part of the group submission:
\begin{displayquote}
    \emph{(Groupwork Quote)} Our group found that proof 1 was the most confusing. It seems as if there needed to be more justifications and explanations, like including that $p/q$ is a reduced fraction. Some of my sentences were not clear, such as my opening.
\end{displayquote}
\vskip.2cm 

In retrospect, it would have been better to ask students to review the original proofs submitted by the group, identify the best traits in them, and then come together as a group to write a stronger version of the proof together. Such asset-based perspective might mitigate risks of embarrassment for emerging proof-writers. \\

\item \emph{The group revision may not demonstrate the group's full potential in writing proofs.}

In the case of the example submission investigated in Section \ref{sec: PE1 success story}, we notice that the group revision of Nora's original proof shows significant improvements with respect to the original proof. We observe, however, that some shortcomings still remain. 

As instructors, we hoped that the revised proof would incorporate all the best features of the four original proofs, but noticed with some disappointment that there were several positive  elements in David's, Yonathan's, and Oliver's original proofs that did not get transferred to the collective revision of Nora's proof. This is not due to a lack of diligence from the students, but rather to the specific language that we (instructors) used in the instructions of the activity.  Specifically, the directions read:  
\vskip.2cm 

\begin{displayquote}
    \emph{(Prompt)} After discussion, identify a part of 1 proof that was confusing. Explain in a few sentences why some group members were confused by this portion. Then, come up with a better written explanation as a group.
\end{displayquote}
\vskip0.2cm

 Students are instructed to focus on fixing \emph{one portion} of a proof that they found confusing, as opposed to joining efforts to write the best collective proof.   
 
 Listed below are some examples of  ``missed opportunities'', namely positive features that are present in one of the  students' original proof but not in the group revision. An extract from the original student work is also included. 
\begin{itemize}
\item The remark that the lemma is applicable because $x$ is prime:
\vskip.2cm 
\begin{displayquote}
    \emph{(Quote from Yonathan)} Since $x$ is prime and $x \mid p\cdot p$, we can use the lemma to state\dots 
\end{displayquote}
    \vskip.2cm

\item The observation that the contradiction follows from the fact that $x>1$, which -- in turn -- follows from the assumption that $x$ is prime. Quoting Yonathan's proof:
\vskip.2cm 
\begin{displayquote}
    \emph{(Quote from Yonathan)}  Since $x$ is prime, $x$ is greater or equal than 2, and  $x$ is a common divisor of $p$ and $q$. However, this is a contradiction because $\frac{p}{q}$ is supposed to be a reduced fraction, meaning gcd$(p,q)=1$.
\end{displayquote}
\vskip.2cm

\item The remark that the assumption that  $\sqrt{x}$ is a reduced fraction can be made ``without loss of generality'':
\vskip.2cm 
\begin{displayquote}
    \emph{(Quote from Oliver)} By contradiction, if $\sqrt{x}$ is rational then 
 $\sqrt{x} =\frac{a}{b}$ where $a$ and $b$ are integers and $b\neq 0$. WLOG , let's assume that $a/b$ is a reduced fraction. 
\end{displayquote}
    \vskip.2cm 

\end{itemize}
\vskip.1cm

In conclusion, while the activity seems to have been useful for students to read and revise their peers' proofs, there is potential for further benefit via a secondary proof review.
One option for future implementations of this activity is for the instructor to give detailed feedback on the group's final proof.
A secondary option is to have groups exchange their finalized proofs with another group and do a collective round of group peer review as a follow-up in-class activity. \\

\item \emph{The assignment is lengthy.}

The intention of this collaborative group review is to have students do some individual thinking about what makes a ``good" proof before meeting with their peers, and see what common attributes appeared. However, the activity spans 7 pages, with only the last 2 submitted for grading. 
This can be daunting for students, and the activity could be improved by providing a list of ``thinking prompts" as a single half-page before the group discussion, thus shortening the physical assignment to a concise 3-4 pages in total. \\

\item \emph{Giving the assignment as a take-home group assignment places additional administrative burden upon students.}

A common theme in student evaluations for this course was some sense of displeasure among students regarding these groupwork take-home assignments. Specifically, students were frustrated by the difficulties of finding a group and of coordinating a convenient meeting time with their group members. 
The class sizes for the intro-to-proofs course are too large for coordination of student groups by the instructor to be feasible outside of assigned lecture and discussion time periods (60+ students).

One implementation suggestion for instructors with similarly large class sizes is to post the individual proof exercise and thinking prompts as a take-home assignment, and to give the collaborative proof-revision and submission portion of the assignment during one of the discussion (recitation) sections of the course.

\end{enumerate}

\section{Proof Evaluation 2: Evaluating Online Resources}

\subsection{Description of Activity}
In Proof Evaluation 2, students are instructed to read a proof of the division algorithm and outline the major steps used. In addition, they are asked to find another proof on the internet that they prefer, and identify what aspects made this other proof better. The main student learning objectives of Proof Evaluation 2 are to:
\begin{enumerate}
    \item Break down complex mathematical arguments.
    \item Evaluate online resources to supplement students' learning.
\end{enumerate}
The tasks in Proof Evaluation 2 are part of a guided process to give students practice with reading, annotating, and then re-structuring the main points of a proof. This is an exercise that every mathematician engages with, both when trying to read and understand the work of others, as well as to structure their own arguments into a series of readable lemmas, theorems, and corollaries. 
By including a reflective activity involving the use of outside resources, the hope is that students will begin to evaluate arguments on the web more critically as they continue their studies. 
Indeed, this practice of evaluating online resources for both clarity and relevance is one that is emphasized deeply in other academic fields, but is often lacking in the undergraduate mathematics curriculum. 
However, this evaluation of online resources is an invaluable skill for professional mathematicians and deserves further attention within upper-division mathematics education.

Initially in this proof evaluation, all student groups are provided with a link to a particular proof of the division algorithm, hosted at Emory University.\footnote{\url{https://mathcenter.oxford.emory.edu/site/math125/proofDivAlgorithm/}} 
The chosen proof is quite long, with many explanatory remarks, and with all relevant cases being treated in equal detail.

Students are then asked to:
\begin{enumerate}[label=(\alph*), ref=(\alph*)]
\item \label{steps in Emory proof} identify the major steps in the proof; 
\item \label{discuss Emory proof} discuss the portions of the proof which are provided to help the reader as opposed to being logically necessary, and assess whether these remarks are helpful or confusing; and 
\item \label{locate internet proofs} locate two additional proofs of the division algorithm, and further discuss one of the found proofs which they collectively agree is ``easy to understand".
\end{enumerate}

The proof hosted at Emory University is correct, as likely were the majority of the proofs found on the internet in part~\ref{locate internet proofs}.  In this sense, Proof Evaluation~2 differs from Proof Evaluation~1 (where the students are assessing each others' proofs in part for correctness).
In the wording of Student Learning Objective~(\ref{SLO correctness and clarity}), this current proof evaluation is more concerned with \emph{clarity} than with \emph{correctness}; though this was originally emphasized in Proof Evaluation 1, we continue this investigation here by asking students what aspects of a proof make it ``easier to understand".

This activity was assigned two weeks after Proof Evaluation 1 was due, in the same format: posted online, and assigned as a take-home assignment for students to complete in one week and submit online. In addition to the material covered before Proof Evaluation 1, students had also covered the basics of functions, divisors and the Euclidean algorithm, and had just started the material on mathematical induction.
To view a copy of the activity as given to the students, please see Appendix \ref{appendix: Proof Eval 2}.

\subsection{Evaluating Online Proofs: When is a proof \emph{easier} to understand?}


Proof Evaluation 2 begins by providing a common baseline example of a proof of the division algorithm, the Emory proof linked above.  By working together as a group to compare the Emory proof with two additional proofs found on the internet, the students establish a context for evaluating further online proofs.

To start the groups on the process of establishing their own ideas of what makes a proof easy to understand, with reference to the Emory proof, part~\ref{discuss Emory proof} instructs the groups to ``identify qualities of this proof that are unnecessary, but are helpful for learning".  This calls on the groups to consider not only whether  the proof is correct versus incorrect, but also less absolute notions, such as whether parts of the proof are \emph{unnecessary but helpful}.

Referring to the attention given to the well-ordering principle in the Emory proof, one group wrote:

\begin{displayquote}
    \emph{(Student Quote)} I found it very helpful, as it helps see the logic of the proof, and that there is small little goals within each proof they are working towards in order to build the bigger picture. I thought it was enlightening, as it gave some insight onto where the proof was headed, instead of just continuously spitting math and logic until it was proven.
\end{displayquote}

\begin{figure}[htp]
    \centering
\includegraphics[width=11cm]{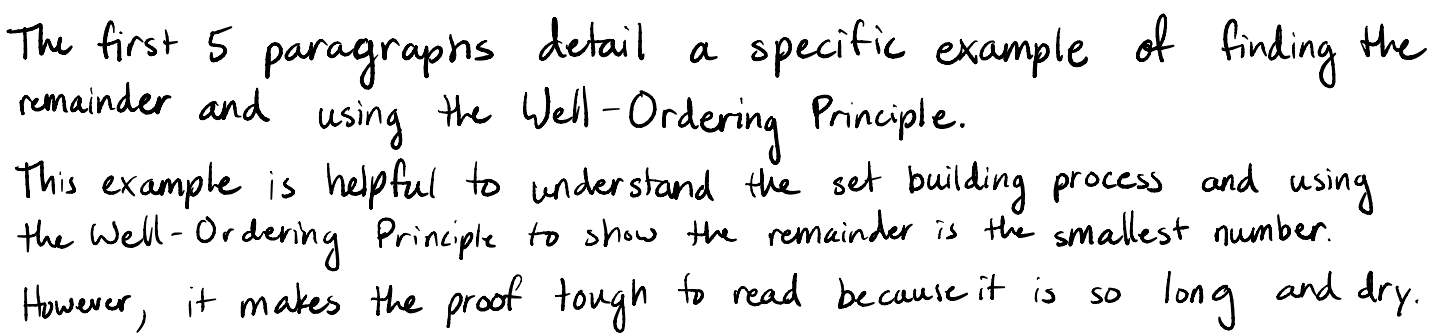}
    \caption{Group submission to part~\ref{discuss Emory proof}}
    \label{fig:5pars}
\end{figure}

Many of the responses to part~\ref{discuss Emory proof} point to a natural tension between providing clarifying comments and  making the proof more overwhelming for the reader in terms of volume.  For example, another group  wrote (see Figure \ref{fig:5pars}),

\begin{displayquote}
    \emph{(Student Quote, Figure \ref{fig:5pars})} The first 5 paragraphs detail a specific example of finding the remainder and using the Well-Ordering Principle.  This example is helpful to understand the set building process and using the Well-Ordering Principle to show the remainder is the smallest number. 
 However, it makes the proof tough to read because it is so long and dry.
\end{displayquote}

In part~\ref{locate internet proofs}, the student groups search the internet for two further proofs of the division algorithm.  They are then asked to single out one of these proofs and explain what makes it ``easier to understand".  (The wording is deliberately vague with regards to whether it is easier to understand than the Emory proof or only easier to understand than the second found proof.)

With specific regards to what makes a proof easier to understand, the different student groups suggested a wide variety of factors.  Here we list the factors identified by the groups of students in our classes, organized by category (and sub-category). The numbers in parentheses indicate the frequency.  
\begin{itemize}
\item \emph{Format/medium/visuals} (11)
\item \emph{Logical contents of the proof / proof type} (6)
    \begin{itemize}
    \item Fewer cases (2)
    \item Contradiction (2)
    \item Proof by cases (1)
    \item Induction (1)
    \end{itemize}
\item \emph{How the proof is written (narrative structure)} (14)
    \begin{itemize}
    \item Clear step-by-step progression (9)
    \item Examples (3)
    \item Initial goal/overview (2)
    \end{itemize}
\item \emph{How the proof is written (terminology)} (13)
    \begin{itemize}
    \item Effective notation (5)
    \item Amount of algebra/arithmetic (4)
    \item Language (4)
    \end{itemize}
\item \emph{How the proof is written (amount of explanation)} (16)
\begin{itemize}
    \item Less to read, more concise (10)
    \item Comprehensive (6)
    \end{itemize}
\end{itemize}

Occasionally a group identified a proof as being incomplete, but beyond that, no mathematical errors were identified by the student groups in any of the resources found online.  The considerations in this proof evaluation assignment are not about correctness, but are instead concerned with less precisely-defined notions, such as whether a sentence in a proof is enlightening or confusing.  These notions are most effectively considered with a group of peers.

The question of what makes a proof easier to understand is subtle.  Student learners may find the Emory proof quite easy to understand, as it covers all steps in complete detail and includes many explanatory comments.  An experienced mathematician, familiar with similar proofs, would likely find much of the Emory proof distracting.  Whether a proof is \emph{easy to understand} and whether it is \emph{clear or confusing} is a social concept, in that it is heavily dependent on the audience.

As an illustration of this, consider the following response, which compares the Emory proof to a proof hosted at Florida State University (FSU).\footnote{\url{https://www.math.fsu.edu/~pkirby/mad2104/SlideShow/s5_1.pdf}(see pages 2-3).}  We will refer to the authors of this response  as Sara's group.

\begin{figure}[htp]
    \centering
\includegraphics[width=11cm]{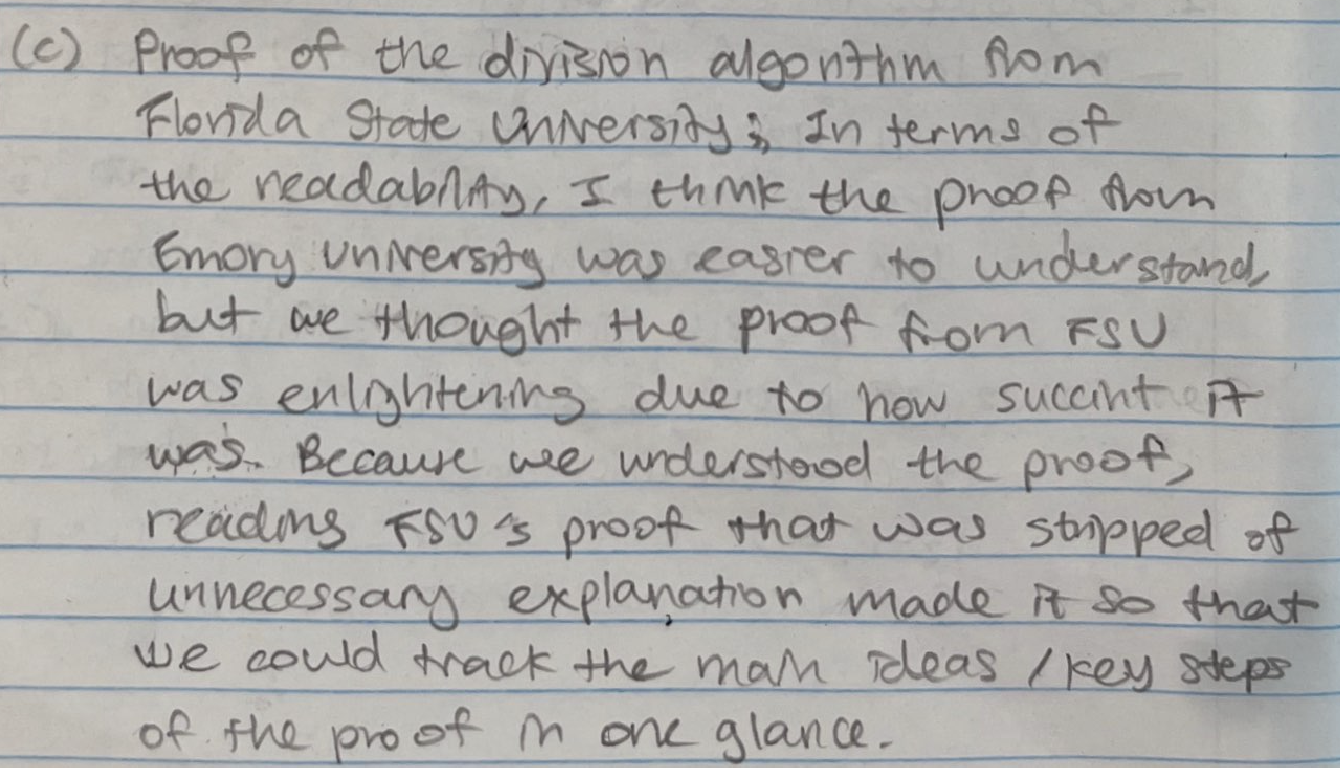}
    \caption{Sara's proof}
    \label{fig:sara}
\end{figure}

Quoting Sara's group's submission: 
\begin{displayquote}
    \emph{(Student Quote, see Figure \ref{fig:sara})} In terms of readability, I think the proof from Emory University was easier to understand, but we thought the proof from FSU was enlightening due to how succinct it was.  Because we understood the proof, reading FSU's proof that was stripped of unnecessary explanation made it so that we could track the main ideas / key steps of the proof in one glance.
\end{displayquote}

We find this response particularly interesting, because the students in Sara's group are seemingly playing \emph{two roles}  within this single assignment.  They are playing the role of student learners, recognizing that the Emory proof was ``easier to understand".  But within the same paragraph, they then take on the role of a more experienced reader.  This role was accessible to Sara's group because they already ``understood the proof".

The students in Sara's group exhibited a very high degree of awareness. They found the Emory University proof easier to understand, but also recognized that the key steps of the proof can be conveyed in a much more concise way.  Their response also implies that their appreciation for the FSU proof is due in part to already understanding the argument.


\subsection{Implementation and Adaptation} 

Discussing the responses as a class would have provided the opportunity to spotlight a number of aspects of the groups' selections.
For instance, numerous groups pointed to a YouTube video from the channel \emph{Mathematical Visual Proofs} \footnote{\url{https://youtu.be/lEcyafR51tM?si=m2DgcHA0I3bs_W_n}}
as an example of a proof which is easier to understand. We feel that  that students generally lack experience   watching video proofs critically, as the choice of  this particular video demonstrates.  For example, one group credits the video with proving the Well-Ordering Principle, but in fact, the video only briefly applies this principle.  Another group recognizes this aspect, and writes: ``While a little more pacing and a deeper dive into well-ordering may be necessary, this is a great video for understanding the bigger picture."  Our point in singling out the use of well-ordering in this particular  video is not to criticize it, nor even to agree that more explanation would be better, but instead to point out that when watching a video proof, as opposed to reading a written proof, it is perhaps more difficult to evaluate the content critically.  (For example, we admit that one of the authors of this article had to watch the video multiple times before understanding how well-ordering was used in the argument.)

Proof Evaluation 2 could  be adapted in other ways to evaluate online materials and proof clarity. For example, students could be given a video ``proof" of some theorem, and then asked to write a formal proof based on the video. There could be a follow-up discussion on which aspects of the video vs.~written proof seem better for understanding a proof. 
Alternatively, students could be presented with one long, detailed proof meant for teaching and one very concise proof, and asked to synthesize the two, keeping only the details they feel are most important.

A possible follow-up to this activity could be to discuss Mochizuki's claimed proof from August 2012 of the \emph{abc conjecture}, which is an instance of an eminent mathematician publishing a claimed proof of a famous open problem.  Mochizuki's proof is famously written in a manner that even experts working in the same area were not able to follow.  Now more than ten years later, Mochizuki's paper has been published (which, barring errors, typically means the paper has been endorsed as being correct), but many (perhaps even the majority of) mathematicians are skeptical of its correctness.  A highly recommended general audience overview of the state of this proof as it stood in 2013 can be found in \cite{mochizuki}.

\section{Concluding Remarks}
Overall, the two activities presented in this paper represent a first attempt to support the development of sociomathematical norms for students in an introduction to proofs course. 
These activities emphasize clarity of argumentation, collaboration, and the evaluation of online resources.
Instructors interested in implementing similar activities can use our Proof Evaluation tasks (provided in the appendix) as a starting point for developing their own activities. 

Notable changes that we suggest are to deliver the activities in-class wherever possible, or to assign the individual portions as take-home exercises and organize the collaborative groupwork as a part of in-person lecture or recitation. Instructors should also consider rephrasing some prompts in a positive, asset-based perspective. Lastly, we suggest performing follow-up activities that emphasize some of the common traits of what students consider to be a ``good" proof, and highlights some of the tension between these traits (e.g., the desire for a preview or outline of the proof vs. the desire for concision). 
\vspace{1cm}


\newpage

\newpage

\printbibliography

\newpage
\appendix
\begin{appendices}

\markboth{Appendix}{Appendix}


{\section{Proof Evaluation 1}

\label{appendix: Proof Eval 1}}

\normalsize
\section*{Instructions}
First complete the \textbf{Proof Exercise} below on your own. Then, exchange proofs with your group members (e.g. by emailing photos of your completed proof). 
Then, complete the \textbf{Review} section on your own, and finally meet up to complete the \textbf{Group Discussion} together.
\vspace{0.2cm}

Each group should fill out and submit a single submission (see end of this document, starting at page 6).

\section*{Proof Exercise}
\textcolor{red}{Prove or disprove: If $x$ is a prime number, then $\sqrt{x}$ is irrational.}
\vspace{0.2cm}

Please make use of the stages of proof-writing we've discussed: interpret, explore, brainstorm, sketch, proof! But send only your completed proof to your group members.
\vspace{0.5cm}

\underline{Write your final proof below:}

\begin{center}
$\updownarrow$\footnote{The symbol $\updownarrow$ stands in for blank space or a page break that was in the original assignment given to students.}
\end{center}

\section*{Review}
\subsection*{Peer Review 1}
Name of Author: 
\vspace{0.3cm}

Name of Reviewer: 
\vspace{0.3cm}

\begin{enumerate}
    \item Do you believe the proof is correct? If not, what about the proof seems incorrect or incomplete?
\[\updownarrow\]
    \item Is the proof satisfying or enlightening? Do you feel that, by reading it, you gain more of an understanding of why the claim is true?
\[\updownarrow\]
    \item Are there any sentences or portions of the proof that were confusing? Why or how were they confusing? (Are there claims that need further justification? Or extra sentences that are not necessary?)
\[\updownarrow\]
\end{enumerate}

\subsection*{Peer Review 2}
Name of Author: 
\vspace{0.3cm}

Name of Reviewer: 
\vspace{0.3cm}

\begin{enumerate}
    \item Do you believe the proof is correct? If not, what about the proof seems incorrect or incomplete?
\[\updownarrow\]
    \item Is the proof satisfying or enlightening? Do you feel that, by reading it, you gain more of an understanding of why the claim is true?
\[\updownarrow\]
    \item Are there any sentences or portions of the proof that were confusing? Why or how were they confusing? (Are there claims that need further justification? Or extra sentences that are not necessary?)
\[\updownarrow\]
\end{enumerate}

\subsection*{Comparison}
Now read your own proof of the claim, and compare it with your two peers’ proofs.

\begin{enumerate}
    \item In what ways are the three proofs similar? Do they all start with the same assumptions? Do they all end with the same conclusion?
\[\updownarrow\]
    \item In what ways are the three proofs different? Do they use different methods of proof? Are some proofs easier to parse than others?
\[\updownarrow\]
    \item Choose one proof that you find easier/clearer than the others. Why do you find it easier to understand? 
\[\updownarrow\]
\end{enumerate}

\section*{Group Discussion}
\begin{enumerate}
    \item First discuss any portions of any proofs that were confusing. Have the author try to convince the group of the validity of their proof/claim in person. 
    \vspace{0.2cm}
    
    After discussion, \textbf{identify a part of 1 proof that was confusing. Explain in a few sentences why some group members were confused by this portion. Then, come up with a better written explanation as a group.}
\[\updownarrow\]

    \item Is there any agreement about which proof was the most clear from your individual 3-way comparisons? Compare the reasons each of you found a particular proof clear/easy to understand. After discussion, \textbf{write a few sentences about what your group agrees upon as qualities of a “good” proof.}
\[\updownarrow\]
\end{enumerate}

\section*{Submission}
First write down the 3 individual proofs your group wrote:
\begin{enumerate}
    \item[] Proof 1:
\[\updownarrow\]
    \item[] Proof 2:
\[\updownarrow\]
    \item[] Proof 3: 
\[\updownarrow\]
\end{enumerate}

Then write down the results of your group discussion by answering the questions below:

\begin{enumerate}
    \item Identify a part of 1 proof that was confusing. Explain in a few sentences why some group members were confused by this portion. Then, come up with a better written explanation as a group. \\
    \textbf{Proof Number:} \\
    \textbf{What was confusing:}
\[\updownarrow\]
    
    \textbf{Revised Proof (Rewrite entire proof, with your group's improvements/corrections):}
\[\updownarrow\]

    \item Write a few sentences about what your group agrees upon as qualities of a “good” proof. Give examples of these qualities from the proofs on the previous page.
\[\updownarrow\]
\end{enumerate}

\newpage

{\section{Proof Evaluation 2}

\label{appendix: Proof Eval 2}}

\normalsize
\section*{Instructions}
In this proof evaluation, we will focus on our ability to \textit{read and understand} mathematics. \\

This group proof evaluation is about the following proof of the division algorithm: \url{https://mathcenter.oxford.emory.edu/site/math125/proofDivAlgorithm/}. \\

Read and understand the proof, then answer the following questions together. \\
\textit{Tip/Recommendation: If you are having a difficult time with digesting the proof, I highly recommend having a pen/paper out and taking small notes while you read. This way, you can also check the smaller steps as you go through the proof. This method of reading proofs really helps me engage with the process; it's very difficult to just read things off of a page and understand them immediately!} \\

Each group should fill out and submit a single submission. Groups can be formed on Canvas under ``People $>$ Groups", and you can join any group number with your group members. Then Canvas will automatically link your group's grades for this submission!

\section*{A Necessary Fact}
For reading and understanding the proof, you will need the following \textbf{axiom}. (Remember that axioms are the underlying assumptions we make, and they do not require proof!).

\begin{axiom}{(Well-Ordering Principle)}
    \emph{} \\
    Every non-empty subset $S$ of the non-negative integers $\Z_{\geq 0}$ has a least element.
\end{axiom}

\begin{example}
\emph{} \\
\begin{itemize}
    \item Let $S = \{2, 4, 6, 8, \ldots\}$ be the set of positive even integers. Then $S$ has a least element, $2$.
    \item Let $S = \{4, 5, 6, 7, 8, 9\}$. Then $S$ has a least element, $4$.
    \item Let $S = \N$, all of the natural numbers. Then $S$ has a least element, $1$.
    \item Let $S = \Z_{\geq 0}$. Then $S$ has a least element, $0$.
\end{itemize}
\end{example}
The well-ordering principle just says we can always choose a unique least element from any set $S \subseteq \Z_{\geq 0}$, as long as $S$ is not empty!
\[\updownarrow\]
\section*{Group Discussion/Submission}
\begin{enumerate}
    \item The proof is quite long. As a group, outline the major steps of the proof: how did the author establish existence, and then uniqueness? Try to identify the key points of each argument.
    \vspace{0.2cm}
    
    \textbf{To establish existence, the key steps were: }
\[\updownarrow\]
    \textbf{To establish uniqueness, the key steps were: }
\[\updownarrow\]
    \item This proof is also written to be used for teaching a math reasoning course like Math 13. As a group, can you identify qualities of this proof that are unnecessary, but are helpful for learning/understanding where the proof idea comes from? \\
    \textbf{Write a few sentences identifying portions of the proof that are not necessary, but might be helpful for a student to read. What makes these portions helpful? When you read the proof, did you find them enlightening or confusing?}
\[\updownarrow\]
    \item Find and read (or watch) at least 2 other proofs of the division algorithm. (E.g. you can Google ``Proofs of the division algorithm", and many things should pop up). As a group, discuss which proofs are easier to understand, and which ones are more confusing. \\
    
    \textbf{Identify one proof your group found easy to understand:} \\
    (If you are writing on paper and want to avoid copying down a long URL, you can give identifying information, for example \textit{proof from Emory, Oxford College, Department of Math and CS}, \\
    or \textit{proof on YouTube, from 3Blue1Brown, with title ``The Division Algorithm"}, etc.)
    \vskip4.5cm

    \textbf{What about this proof made it easier to understand?}
\end{enumerate}
\[\updownarrow\]

\end{appendices}

\end{document}